\documentclass{article}
\usepackage{amsmath,amsthm,amsfonts,eucal}

\numberwithin{equation}{section}

\def\ca{{\mathcal A}}

\def\ck{{\mathcal K}}

\def\cu{{\mathcal U}}

\def\bg{{\mathbb G}}

\def\bn{{\mathbb N}}
\def\br{{\mathbb R}}
\def\bz{{\mathbb Z}}

\def\a{\alpha}
\def\b{\beta}
\def\g{\gamma}        
\def\d{\delta}        \def\D{\Delta}
\def\eps{\varepsilon}

\def\l{\lambda}       
\def\m{\mu}
\def\n{\nu}

\def\t{\tau}

\def\f{\varphi}

\def\om{\omega}        \def\O{\Omega}

\def\e#1{{\rm e}^{#1}}
\def\itm#1{\item[$(#1)$]}

\def\limr{\lim_{r\to\infty}}
\def\limR{\lim_{R\to\infty}}
\def\lsup{\limsup_{R\to\infty}}

\DeclareMathOperator{\inj}{inj}

\DeclareMathOperator{\length}{length}

\def\ad#1{d_\infty(#1)}
\def\md#1{d_0(#1)}

\def\dE{{\partial E}}

\DeclareMathOperator{\Lim}{Lim}

\def\bg{C$^{\infty}$-bounded geometry }

\newtheorem{Thm}{Theorem}[section]
\newtheorem{Cor}[Thm]{Corollary}
\newtheorem{Prop}[Thm]{Proposition}
\newtheorem{Lemma}[Thm]{Lemma}
\theoremstyle{definition}
\newtheorem{Dfn}[Thm]{Definition}
\newtheorem{exmp}[Thm]{Example}
\newtheorem{Assump}[Thm]{Assumption}
\theoremstyle{remark}
\newtheorem{rem}[Thm]{Remark} 
\newtheorem{ack}{Acknowledgement} 
 \title{\huge An asymptotic 
 dimension for metric spaces, and the $0$-th Novikov-Shubin invariant}
 \author{Daniele Guido$^{1}$, Tommaso Isola$^{2}$\\
 $(1)$ Dipartimento di Matematica,\\ Universit\`a della Basilicata,\\ 
 I--85100 Potenza, Italy.\\
 $(2)$ Dipartimento di Matematica,\\ Universit\`a di Roma ``Tor Vergata'',\\ 
 I--00133 Roma, Italy.\\
 {\tt guido@unibas.it,\ isola@mat.uniroma2.it}}
\date{March 24, 2000}
\begin{document}
\maketitle
\markboth{Asymptotic dimension package}
{Asymptotic dimension and Novikov-Shubin invariants}
\renewcommand{\sectionmark}[1]{}
\bigskip

\begin{abstract}
	A nonnegative number $d_{\infty}$, called asymptotic dimension, is 
	associated with any metric space.  Such number detects the 
	asymptotic properties of the space (being zero on bounded metric 
	spaces), fulfills the properties of a dimension, and is invariant 
	under rough isometries.  It is then shown that for a class of open 
	manifolds with bounded geometry the asymptotic dimension coincides 
	with the $0$-th Novikov-Shubin number $\a_{0}$ defined in a 
	previous paper [D. Guido, T. Isola, J. Funct.  Analysis, {\bf 176} 
	(2000)].  Thus the dimensional interpretation of $\a_{0}$ 
	given in the mentioned paper in the framework of noncommutative 
	geometry is established on metrics grounds.  Since the asymptotic 
	dimension of a covering manifold coincides with the polynomial 
	growth of its covering group, the stated equality generalises to 
	open manifolds a result by Varopoulos.
\end{abstract}

\newpage
\setcounter{section}{-1}
\section{Introduction.}\label{sec:intro}

In a recent paper \cite{GI4}, we extended the notion of Novikov-Shubin 
numbers to amenable open manifolds and showed that they have a 
dimensional interpretation in the framework of noncommutative 
geometry.  Here we introduce an asymptotic dimension for metric 
spaces, which is an asymptotic counterpart of the Kolmogorov dimension 
\cite{KT}, and show that for a class of open manifolds it coincides 
with the $0$-th Novikov-Shubin number.

The dimension introduced by Kolmogorov and Tihomirov, also called box
dimension, ``corresponds to the possibility of characterizing the
``massiveness'' of sets in metric spaces by the help of the order of
growth of the number of elements of their most economical
$\eps$-coverings, as $\eps\to0$'' \cite{KT}.  For non totally bounded
sets, denoting by $n(r,R)$ the minimum number of balls of radius $r$
necessary to cover a ball of radius $R$ (and given center), the box
dimension is the ``order of infinite'' of $n(r,R)$ when $r\to0$ (for
big $R$, and often independently of $R$).  

We then define the asymptotic dimension $d_\infty$ as the ``order of
infinite'' of $n(r,R)$ when $R\to\infty$ (for big $r$, and often
independently of $r$),
$$
\ad{X} = \limr\limR \frac{\log n(r,R)}{\log R},
$$
show that $d_\infty$ is a dimension, namely
$d_\infty(X\cup Y)=\max(d_\infty(X),d_\infty(Y))$ and $d_\infty(X\times
Y)\leq d_\infty(X)+d_\infty(Y)$, and prove that it is invariant under
rough isometries.

Finally we show that the asymptotic dimension of an open manifold with 
C$^\infty$-bounded geometry and satisfying an isoperimetric inequality 
introduced by Grigor'yan \cite{Grigoryan94} coincides with the $0$-th 
Novikov-Shubin number $\a_0$ as defined in \cite{GI4}.  On the one 
hand this strengthens the dimensional interpretation given in 
\cite{GI4}, and on the other it shows that the limit procedure used in 
the definition of $\a_0$ does not affect the result.  Moreover, the 
quasi-isometry invariance of the $\alpha_{p}$ proved in \cite{GI4} 
becomes rough isometry invariance for the case of $\a_{0}$.

Since the asymptotic dimension of a manifold with C$^\infty$-bounded geometry 
may be computed in terms of its volume growth, the equality between 
$\a_{0}$ and $d_{\infty}$ may be seen 
as a generalization of the result of Varopoulos \cite{Varopoulos1} for 
covering manifolds, namely the equality between $\a_0$ and the growth 
of the covering group.

We recall that the Novikov-Shubin numbers \cite{NS} where introduced, 
after the definition by Atiyah \cite{Atiyah} of the $L^2$-Betti 
numbers in terms of the von~Neumann trace of the covering group, as 
finer invariants of the spectral behaviour of the $p$-Laplacian near 
zero, and where shown to be homotopy invariant by Gromov and Shubin 
\cite{GS}.  It was observed by Roe \cite{Roe1} that, when the covering 
group is amenable, the von~Neumann trace of an operator may be 
computed as an average of its integral kernel on the manifold w.r.t. a 
suitable exhaustion.  Hence this procedure may be extended to 
manifolds admitting an amenable exhaustion.  We show that, for 
manifolds satisfying Grigor'yan isoperimetric inequality, an amenable 
exhaustion exists and is given by concentric balls of increasing 
radius.

In \cite{GI3} we showed that, from the operator algebraic point of 
view, the step from amenable coverings to amenable exhaustions 
corresponds to passing from a normal semifinite trace on a von~Neumann 
algebra to a semicontinuous semifinite trace on a C$^*$-algebra.  The 
latter does not necessarily contain spectral projections, however, the 
spectral density function may be still defined, making use of the 
noncommutative Riemann integration developed in \cite{GI4}, hence 
Novikov-Shubin numbers are defined.  Moreover, this definition 
coincides, at least for $\a_{0}$, with the definition given in terms 
of the trace of the heat kernel, which does not require Riemann 
integration, and which is used here.  More precisely, 
$$
\a_{0}(M) := 2 \limsup_{t\to\infty}\frac{\t(\e{-t\D_{p}})}{\log 1/t},
$$
where $\t$ is the (semicontinuous semifinite) trace associated with an 
amenable exhaustion of $M$.

We now recall the dimensional interpretation of the Novikov-Shubin 
numbers given in \cite{GI4} in the framework of noncommutative 
geometry.

On the one hand these numbers are defined in terms of the low 
frequency behaviour of the $p$-Laplacians, or the large time behaviour 
of the $p$-heat kernel, therefore they are the large scale counterpart 
of the spectral dimension, namely of the dimension as it is recovered 
by the Weyl asymptotics.

On the other hand, recall that in Alan Connes' noncommutative geometry 
\cite{Co}, a non-trivial singular trace, associated to some power of 
the resolvent of the Dirac operator, plays the role of integration 
over the noncommutative space, and that power of the resolvent is the 
metric dimension of the space. 

This is analogous of what happens in geometric measure theory where a 
dimension is the unique exponent of the diameter of a ball which gives, 
via Hausdorff procedure, a (possibly) non trivial 
measure on the space.

Therefore it is in this context that Novikov-Shubin numbers are 
interpreted as asymptotic spectral dimensions, since it was shown in 
\cite{GI4} that the operator $\Delta_p^{-1/2}$, raised to the power 
$\alpha_p$, is singularly traceable.

The identification of $\a_{0}$ with $d_{\infty}$ proved here puts on 
metric grounds the dimensional character of the $0$-th Novikov-Shubin 
number.

Finally we study the relation of $d_{\infty}$ with the notion of 
(metric) asymptotic dimension for cylindrical ends given by Davies 
\cite{Davies}.  Such definition is given in terms of the volume growth 
of the end, therefore when the end has bounded geometry Davies' 
asymptotic dimension coincides with ours.  Indeed Davies requires the 
growth to be exactly polynomial, therefore, in contrast with ours, his 
dimension is not always defined.  Davies also introduced a class of 
cylindrical ends called standard ends.  We show that for standard ends 
with Davies' asymptotic dimension $D$ the equality $d_{\infty}=D$ 
holds with no further assumptions.  Making use of Davies standard ends 
one observes that $d_{\infty}$ for a open manifold may take any value 
in $[1,\infty]$.  Also, we discuss examples of standard ends where the 
growth is not exactly polynomial.

Some of the results contained in the present paper have been announced 
in several international conferences.  In particular we would like to 
thank the Erwin Schr\" odinger Institute in Vienna, where the first 
draft of this paper was completed, and the organisers of the 
``Spectral Geometry Program'' for their kind invitation.

 \section{An asymptotic dimension for metric spaces.}\label{sec:first}
 
 The main purpose of this section is to introduce an asymptotic 
 dimension for metric spaces.  Other notions of asymptotic dimension 
 have been considered by Gromov \cite{Gromov} (see also the papers by 
 Yu \cite{Yu} and Dranishnikov \cite{Dra}).  Davies \cite{Davies} 
 proposed a definition in the case of cylindrical ends of a Riemannian 
 manifold.  We shall give a definition of asymptotic dimension in the 
 setting of metric dimension theory, based on the (local) Kolmogorov 
 dimension \cite{KT} and state its main properties.  We compare our 
 definition with Davies' in the next Section, and discuss its 
 relations with Gromov's in Remark \ref{gromov}.

 In the following, unless otherwise specified, $(X,\d)$ will denote a 
 metric space, $B_X(x,R)$ the open ball in $X$ with centre $x$ and 
 radius $R$, $n_r(\O)$ the least number of open balls of radius $r$ 
 which cover $\O\subset X$, and $\n_r(\O)$ the largest number of 
 disjoint open balls of radius $r$ centered in $\O$.
 
 Kolmogorov and Tihomirov \cite{KT} defined a dimension for totally 
 bounded metric spaces $X$ as $\md{X}:=\limsup_{r\to0} \frac{\log 
 n_r(X)}{\log(1/r)}$.  A natural extension to all metric spaces is 
 given by $\md{X}=\lim_{R\to\infty}\limsup_{r\to0} \frac{\log 
 n_r(B_X(x,R))}{\log(1/r)}$.  It can be shown that $d_{0}$ is 
 independent of $x$, is a dimension, namely it satisfies the 
 properties of Theorem \ref{Thm:adim}, and is invariant under 
 bi-Lipschitz mappings. 
 
 The following definition gives a natural asymptotic counterpart of 
 the  dimension of Kolmogorov-Tihomirov.

 \begin{Dfn}\label{1.2.1}  
	 Let $(X,\d)$ be a metric space.  We call
	 $$
	 \ad{X}:=\limr\lsup \frac{\log n_r(B_X(x,R)) }{ \log R},
	 $$
	 the {\it asymptotic dimension} of $X$. 
 \end{Dfn}

 Let us remark that, as $n_r(B_X(x,R))$ is nonincreasing in $r$, the 
 function
 $$
 r\mapsto\lsup\frac{\log n_r(B_X(x,R)) }{ \log R}
 $$ 
 is nonincreasing too, so the $\lim_{r\to\infty}$ exists.
 
 \begin{Prop}\label{1.2.3} 
	 $\ad{X}$ does not depend on $x$.
 \end{Prop} 
 \begin{proof}
	 Let $x,y\in X$, and set $\d:=\d(x,y)$, so that $B(x,R)\subset 
	 B(y,R+\d) \subset B(x,R+2\d)$.  This implies,
	 \begin{align*}
		 \frac{\log n_r(B(x,R)) }{ \log R} &\leq \frac{\log 
		 n_r(B(y,R+\d)) }{ \log (R+\d)}\ \frac{\log (R+\d) }{ \log R} 
		 \\
		 &\leq \frac{\log n_r(B(x,R+2\d)) }{ \log (R+2\d)}\ \frac{\log 
		 (R+2\d) }{ \log R}
	 \end{align*}
	 so that, taking $\lsup$ and then $\limr$ we get the thesis.  
 \end{proof}
 
 The following lemma is proved in \cite{KT}.  For the sake of 
 completeness, we include a proof.

 \begin{Lemma}\label{1.1.1} 
	 $n_r(\O) \geq \n_r(\O)\geq n_{2r}(\O)$. 
 \end{Lemma}
 \begin{proof} 
	 For the first inequality, let $B_{X}(x_{i},r)$, 
	 $i=1,\ldots,\n_{r}(\O)$, be disjoint balls with centres in $\O$.  
	 Then any $r$-ball of a covering of $\O$ may contain at most one 
	 of the $x_{i}$'s.  Indeed, if $B_{X}(x,r)\ni x_{i},x_{j}$, then 
	 $B_{X}(x_{i},r)\cap B_{X}(x_{j},r) \supset \{x\} 
	 \neq\emptyset$, so that $x_{i}=x_{j}$.\\
	 As for the second inequality, we need to prove it only when $\n_r$ is 
	 finite.  Let us assume that $\{B(x_i,r)\}_{i=1}^{\n_r(\O)}$ are 
	 disjoint balls centered in $\O$ and observe that, for any $y\in 
	 \O$, $\d(y,\bigcup_{i=1}^{\n_r(\O)} B(x_i,r))<r$, otherwise 
	 $B(y,r)$ would be disjoint from $\bigcup_{i=1}^{\n_r(\O)} 
	 B(x_i,r)$, contradicting the maximality of $\n_r$.  So for all 
	 $y\in \O$ there is $j$ s.t. $\d(y, B(x_j,r))<r$, that is 
	 $\O\subset \bigcup_{i=1}^{\n_r(\O)} B(x_i,2r)$, which implies the 
	 thesis.
 \end{proof}

 \begin{Prop}\label{1.2.4}
	 $$
	 \ad{X}=\limr\lsup \frac{\log \n_r(B_X(x,R)) }{ \log R}
	 $$ 
 \end{Prop}
 \begin{proof} 
	 Follows easily from Lemma \ref{1.1.1}.
 \end{proof}

 \begin{Dfn}\label{1.2.8} 
	 Let $X,Y$ be metric spaces, $f:X\to Y$ is said to be a rough 
	 isometry if there are $a\geq1$, $b,\eps\geq0$ s.t. 
	 \itm{i} $a^{-1}\d_X(x_1,x_2)-b \leq \d_Y(f(x_1),f(x_2)) \leq a 
	 \d_X(x_1,x_2)+b$, for all $x_1,x_2\in X$, 
	 \itm{ii} $\bigcup_{x\in X} B_Y(f(x),\eps) = Y$ 
 \end{Dfn}

 \begin{Lemma}\label{1.2.9} {\rm (\cite{Chavel}, Proposition 4.3)} 
	 If $f:X\to Y$ is a rough isometry, there is a rough isometry 
	 $f^-:Y\to X$, with constants $a,b^-,\eps^-$, s.t. 
	 \itm{i} $\d_X(f^-\circ f(x),x)<c_X$, $x\in X$, 
	 \itm{ii} $\d_Y(f\circ f^-(y),y)<c_Y$, $y\in Y$. 
 \end{Lemma}

 \begin{Thm}\label{1.2.10} 
	 Let $X,Y$ be metric spaces, and $f:X\to Y$ a rough isometry.  
	 Then $\ad{X}=\ad{Y}$.
 \end{Thm}
 \begin{proof}
	 Let $x_0\in X$, then for all $x\in B_X(x_0,r)$ we have
	 $
	 \d_Y(f(x),f(x_0))\leq a \d_X(x,x_0)+b<ar+b
	 $ 
	 so that $f(B_X(x_0,r))\subset B_Y(f(x_0),ar+b)$. 
	 Then, with $n:=n_r(B_Y(f(x_0),aR+b))$, 
	 $$
	 f(B_X(x_0,R)) \subset \bigcup_{j=1}^{n} B_Y(y_j,r),
	 $$
	 which implies 
	 \begin{align*}
		 f^-\circ f(B_X(x_0,R)) & \subset \bigcup_{j=1}^{n} 
		 f^-(B_Y(y_j,r)) \\
		 & \subset \bigcup_{j=1}^{n} B_X(f^-(y_j),ar+b^-).
	 \end{align*}
	 Let $x\in B_X(x_0,R)$, and $j$ be s.t. $f^-\circ f(x)\in 
	 B_X(f^-(y_j),ar+b^-)$, then
	 $$
	 \d_X(x,f^-(y_j))\leq \d_X(x,f^-\circ f(x))+\d_X(f^-\circ 
	 f(x),f^-(y_j))<c_X+ar+b^-,
	 $$ 
	 so that  
	 $$
	 B_X(x_0,R)\subset \bigcup_{j=1}^{n} B_X(f^-(y_j),ar+b^-+c_X),
	 $$
	 which implies $n_{ar+b^-+c_X}(B_X(x_0,R))\leq 
	 n_r(B_Y(f(x_0),aR+b))$.  \\
	 Finally 
	 \begin{align*}
		 \ad{X} & = \limr\lsup\frac {\log n_r(B_X(x_0,R))}{\log R} \\
		 & = \limr\lsup\frac {\log n_{ar+b^-+c_X}(B_X(x_0,R))}{\log R}\\
		 & \leq \limr\lsup\frac {\log n_r(B_Y(f(x_0),aR+b))}{\log R} \\
		 & = \limr\lsup\frac {\log n_r(B_Y(f(x_0),R))}{\log R} \\
		 & = \ad{Y}
	 \end{align*}
	 and exchanging the roles of $X$ and $Y$ we get the thesis.
 \end{proof}

 \begin{Thm}\label{Thm:adim} 
	 The set function $d_\infty$ is a dimension, namely it satisfies 
	 \item{$(i)$} If $X\subset Y$ then $\ad{X}\leq \ad{Y}$.  
	 \item{$(ii)$} If $X_1,X_2\subset X$ then $\ad{X_1\cup X_2} = \max 
	 \{ \ad{X_1}, \ad{X_2} \}$.  
	 \item{$(iii)$} If $X$ and $Y$ are metric spaces, then 
	 $\ad{X\times Y} \leq \ad{X}+\ad{Y}$.
 \end{Thm}
 \begin{proof}
	 $(i)$ Let $x\in X$, then $B_X(x,R)\subset B_Y(x,R)$ and the claim 
	 follows easily.  \\
	 $(ii)$ By part $(i)$, we get $\ad{X_1\cup X_2} \geq \max \{ 
	 \ad{X_1},\ad{X_2} \}$.  So we need to prove the converse 
	 inequality, only when both $X_{1}$ and $X_{2}$ have finite 
	 asymptotic dimension.  Let $x_i\in X_i$, $i=1,2$, and set 
	 $\d=\d(x_1,x_2)$, $a=\ad{X_1}$, $b=\ad{X_2}$, with e.g. $a\leq 
	 b<\infty$.  Then $B_{X_1\cup X_2}(x_1,R) \subset B_{X_1}(x_1,R) 
	 \cup B_{X_2}(x_2,R+\d)$, therefore
	 \begin{equation}\label{ineq1}
		 n_r( B_{X_1\cup X_2}(x_1,R) ) \leq n_r( B_{X_1}(x_1,R) ) +
		 n_r( B_{X_2}(x_2,R+\d) ).
	 \end{equation}
	 Besides, $\forall \eps>0$  $\exists r_{0}>0$ s.t. $\forall 
	 r>r_{0}$ $\exists R_0=R_0(\eps,r)$ s.t. $\forall R>R_0$
	 \begin{align*}
		 n_r(B_{X_1}(x_1,R))  & \leq R^{a+\eps} \\
		 n_r(B_{X_2}(x_2,R+\d)) & \leq R^{b+\eps},
	 \end{align*}
	 hence, by inequality (\ref{ineq1}),
	 $$
	 n_r(B_{X_1\cup X_2}(x_1,R)) \leq R^{a+\eps}+R^{b+\eps} = 
	 R^{b+\eps}(1+R^{a-b}).
	 $$
	 Finally,
	 $$
	 \frac {\log n_r(B_{X_1\cup X_2}(x_1,R))} {\log R} \leq b+\eps+ 
	 \frac {\log(1+R^{a-b})}{\log R}.
	 $$
	 Taking the $\lsup$ and then the $\limr$ we get
	 $$
	 \ad{X_1\cup X_2}\leq \max\{\ad{X_1},\ad{X_2}\}+\eps
	 $$
	 and the thesis follows by the arbitrariness of $\eps$.  \\
	 $(iii)$ By  Proposition \ref{1.2.10}, we may endow $X\times Y$ with any 
	 metric roughly isometric to the product metric, e.g.
	 \begin{equation}\label{e:prodmetric}
	 	\d_{X\times Y}((x_1,y_1), (x_2,y_2))=\max\{ 
	 \d_X(x_1,x_2),\d_Y(y_1,y_2)\}.
	 \end{equation}
	 Then, by $n_r(B_{X\times Y}((x,y),R))\leq n_r(B_X(x,R))\ 
	 n_r(B_Y(y,R))$, the thesis follows easily.  
 \end{proof}

 \begin{rem}\label{1.2.7} 
	 $(a)$ In part $(ii)$ of the previous theorem we considered $X_1$ and 
	 $X_2$ as metric subspaces of $X$.  If $X$ is a Riemannian 
	 manifold and we endow the submanifolds $X_1$, $X_2$ with their 
	 geodesic metrics this property does not hold in general.  A 
	 simple example is the following.  Let $f(t) := (t\cos t,t\sin 
	 t)$, $g(t) := (-t\cos t, -t\sin t)$, $t\geq0$ planar curves, and 
	 set $X,\ Y$ for the closure in $\br^2$ of the two connected 
	 components of $\br^2\setminus (G_f\cup G_g)$, where $G_f,\ G_g$ 
	 are the graphs of $f,\ g$, and endow $X,\ Y$ with the geodesic 
	 metric.  Then $X$ and $Y$ are roughly-isometric to $[0,\infty)$ 
	 (see below) so that $\ad{X}=\ad{Y}=1$, while $\ad{X\cup Y}=2$. \\
	 $(b)$ The choice of the $\limsup$ in 
	 Definition~\ref{1.2.1} is the only one compatible with the 
	 classical dimensional inequality stated in Theorem~\ref{Thm:adim} 
	 $(iii)$.  
 \end{rem}

 In what follows we show that when $X$ is equipped with a suitable 
 measure, the asymptotic dimension may be recovered in terms of the 
 volume asymptotics for balls of increasing radius.  This is analogous 
 to the fact that the local dimension may be recovered in terms of the 
 volume asymptotics for balls of infinitesimal radius.

 \begin{Dfn}\label{1.2.12} 
	 A Borel measure $\m$ on $(X,\d)$ is said to be uniformly bounded 
	 if there are functions $\b_1,\b_2$, s.t. $0<\b_1(r)\leq 
	 \m(B(x,r)) \leq \b_2(r)$, for all $x\in X$, $r>0$.  \\
	 That is $\b_1(r):= \inf_{x\in X} \m(B(x,r)) >0$, and $\b_2(r) := 
	 \sup_{x\in X} \m(B(x,r)) <\infty$.
 \end{Dfn}

 \begin{Prop}\label{1.2.13} 
	 If $(X,\d)$ has a uniformly bounded measure, then every ball in 
	 $X$ is totally bounded (so that if $X$ is complete it is locally 
	 compact).
 \end{Prop}
 \begin{proof}
	 Indeed, if there is a ball $B=B(x,R)$ which is not totally 
	 bounded, then there is $r>0$ s.t. every $r$-net in $B$ is 
	 infinite, so $n_r(B)$ is infinite, and $\n_r(B)$ is infinite too.  
	 So that $\b_2(R)\geq \m(B) \geq \sum_{i=1}^{\n_r(B)} \m(B(x_i,r)) 
	 \geq \b_1(r)\n_r(B) = \infty$, which is absurd.  
 \end{proof}

 \begin{Prop}\label{1.2.14} 
	 If $\m$ is a uniformly bounded Borel measure on $X$ then
	 $$
	 \ad{X}=\lsup\frac{\log\m(B(x,R))}{ \log R} .
	 $$ 
 \end{Prop}
 \begin{proof}
	 As $\bigcup_{i=1}^{\n_r(B(x,R))} B(x_i,r) \subset B(x,R+r) 
	 \subset \bigcup_{j=1}^{n_r(B(x,R+r))} B(y_j,r)$, we get
	 $$
	 \b_2(r)n_r(B(x,R+r))\geq \m(B(x,R+r)) \geq \b_1(r)\n_r(B(x,R)) 
	 \geq \b_1(r)n_{2r}(B(x,R)),
	 $$
	 by Lemma \ref{1.1.1}. So that 
	 $$
	 \b_1(r/2)\leq \frac{\m(B(x,R+r/2))}{ n_r(B(x,R))}, \qquad 
	 \frac{\m(B(x,R))}{ n_r(B(x,R))} \leq \b_2(r),
	 $$
	 and the thesis follows easily.
 \end{proof}

 Let us now consider the particular case of complete Riemannian 
 manifolds.
 
  \begin{Dfn}\label{2.1.1}
	 Let $(M,g)$ be an $n$-dimensional complete Riemannian manifold.  
	 We say that $M$ has bounded geometry if it has positive 
	 injectivity radius, sectional curvature bounded from above, and 
	 Ricci curvature bounded from below.
 \end{Dfn}
 
 \begin{Lemma}\label{2.1.2} 
	 Let $M$ be an $n$-dimensional complete Riemannian manifold with 
	 bounded geometry.  Then the Riemannian volume is a uniformly 
	 bounded measure.
 \end{Lemma}
 \begin{proof}
	 We can assume, without loss of generality, that the sectional 
	 curvature is bounded from above by some positive constant $c_1$ 
	 and the Ricci curvature is bounded from below by $(n-1)c_2 g$, with 
	 $c_2<0$.  Then, denoting with $V_\d(r)$ the volume of a ball of 
	 radius $r$ in a manifold of constant sectional curvature equal to 
	 $\d$, we can set $\b_1(r) := V_{c_1}(\min\{r, r_0\})$, and $\b_2:= 
	 V_{c_2}(r)$, where $r_0:= \min\{ \inj(M), \frac{\pi}{\sqrt{c_1}} 
	 \}$, and $\inj(M)$ is the injectivity radius of $M$.  Then the 
	 result follows from (\cite{Chavel}, p.119,123).
 \end{proof}

 \begin{Prop}\label{2.1.3}
	 Let $M,N$ be complete Riemannian manifolds.  
	 \itm{i} If $M$ is non-compact, then $\ad{M}\geq 1$
	 \itm{ii} If $M$ has bounded geometry, then 
	 $$
	 \ad{M} = \lsup
	 \frac{\log vol(B_{M}(x,R))}{\log R},\ x\in M
	 $$
	 \itm{iii} If $M,N$ have bounded geometry, and
	 admit asymptotic dimension in a strong sense, that is
	 $\ad{M} = \lim_{R\to\infty} \frac{\log vol(B_{M}(x,R))}{\log R}$,  
	 $x\in M$, and analogously for $N$, then 
	 $$
	 \ad{M\times N}=\ad{M}+\ad{N}.
	 $$
 \end{Prop}
 \begin{proof}
	 $(i)$ Let us fix $x_{0}\in M$ and $R>0$, and consider $x_{R}\in 
	 M$ s.t. $\d(x_{0},x_{R})=R$, which exists because $M$ is not 
	 compact, and let $\g:[0,1]\to M$ be a minimizing geodesics 
	 between $\g(0)=x_{0}$, $\g(1) = x_{R}$.  Clearly 
	 $\g([0,1))\subset B_{M}(x_{0},R)$, hence, if $x_{1},\ldots,x_{k}$ 
	 are the centres of a minimal covering by $r$-balls of 
	 $\g([0,1))$, we have $k\leq n_{r}( B_{M}(x_{0},R) )$. Then 
	 $$
	 R= \length(\g) \leq \sum_{i=1}^{k}\length(\g\cap 
	 B_{M}(x_{i},r))\leq 2rk,
	 $$
	 namely $ n_{r}(B_{M}(x_{0},R) ) \geq \frac{R}{2r}$, from which 
	 the thesis follows.  \\
	 $(ii)$  The result follows from Lemma \ref{2.1.2} and 
	 Proposition~\ref{1.2.14}.  \\
	 $(iii)$ As in the proof of Theorem \ref{Thm:adim} $(iii)$, we may 
	 endow $M\times N$ with the metric (\ref{e:prodmetric}).  Then 
	 $vol(B_{M\times N}((x,y),R))=vol(B_M(x,R)) vol(B_N(y,R))$, and we 
	 get
	 \begin{align*}
		 \ad{M\times N} &= \limR \frac{\log vol (B_{M\times 
		 N}((x,y),R))}{ \log R} \\
		 & = \limR \frac{\log vol (B_M(x,R))}{ \log R} + \limR 
		 \frac{\log vol (B_N(y,R))}{ \log R} \\
		 & = \ad{M}+\ad{N}.
	 \end{align*}
 \end{proof}

 \begin{rem}
	 $(a)$ Conditions under which the inequality in Theorem~\ref{Thm:adim} 
	 $(iii)$ becomes an equality are often studied in the case of 
	 (local) dimension theory (cf.  \cite{Pontriagin,Salli}).  The 
	 previous Proposition gives such a condition for the asymptotic 
	 dimension. \\
	 $(b)$ As the asymptotic dimension is invariant under rough 
	 isometries, it is natural to substitute the continuous space with 
	 a coarse graining, which destroys the local structure, but 
	 preserves the large scale structure.  Then (cf.  \cite{Chavel}, 
	 Theorem 4.9) if $M$ is a complete Riemannian manifold with Ricci 
	 curvature bounded from below, $M$ is roughly isometric to any of 
	 its discretizations, endowed with the combinatorial metric.  
	 Therefore $M$ has the same asymptotic dimension of any of its 
	 discretizations.  The previous result, together with the 
	 invariance of the asymptotic dimension under rough isometries, 
	 shows that, when $M$ has a discrete group of isometries $\Gamma$ 
	 with a compact quotient, the asymptotic dimension of the manifold 
	 coincides with the asymptotic dimension of the group, hence with 
	 its growth (cf.  \cite{GI5}).  Therefore, by a result of 
	 Varopoulos \cite{Varopoulos1}, it coincides with the 0-th 
	 Novikov-Shubin invariant.  We will generalise this result in 
	 Section \ref{subsec:NSinvariant}.
 \end{rem}

 Let us conclude this Section with some examples. Other examples are 
 contained in the next Section.

 \begin{exmp}\label{1.2.16} 
	 \itm{i} $\br^n$ has asymptotic dimension $n$. 
	 \itm{ii} Set $X:= \cup_{n\in\bz}\{(x,y)\in\br^2 : 
	 \d((x,y),(n,0))<\frac{1}{4} \}$, endowed with the Euclidean 
	 metric, then $\md{X}=2$, $\ad{X}=1$.  
	 \itm{iii} Set $X=\bz$ with the counting measure, then $\md{X}=0$, 
	 and $\ad{X}=1$.  
	 \itm{iv} Let $X$ be the unit ball in an infinite dimensional 
	 Banach space.  Then $\md{X}=+\infty$ while $\ad{X}=0$.
	 \itm{v} Let $X$ be the $\bz^{\infty}$-lattice determined by an 
	 orthonormal base in  an infinite dimensional 
	 Hilbert space. Then $\md{X} =0$ while $\ad{X} = \infty$.
 \end{exmp}

 \begin{rem}\label{gromov} M. Gromov introduced a notion of "large 
 scale dimension" for metric spaces: the asymptotic dimension of $X$ 
 is the smallest integer $n$ such that, for any $r>0$, there is a 
 cover ${\cu} = \{U_i\}$ of $X$ such that the diameters of the sets 
 $U_i$ are bounded, and no ball of radius $r$ meets more than $n+1$ of 
 them.
 
 Our asymptotic dimension can be very different from Gromov's.  For 
 example hyperbolic space $H_{n}$ has finite Gromov dimension, but 
 $\ad{H_{n}}=\infty$.  Conversely, one can find a sequence of 
 cylindrical ends with fixed $d_{\infty}$ and arbitrarily large Gromov 
 dimension (cf. Corollary \ref{2.3.4}).  
 
 The two notions however, coincide on some very special spaces, such 
 as cartesian products of $\br^{n}$ and a compact set with the product 
 metric.  Moreover both dimensions are in a sense ``coarse'', since 
 they are invariant under rough isometries.
 
 Finally we remark that Gromov dimension is an asymptotic topological 
 dimension, since it is a coarse analogue of the Lebesgue covering 
 dimension, according to Dranishnikov \cite{Dra}.  Ours instead is an 
 asymptotic metric dimension.  Indeed it is an asymptotic counterpart 
 of the Kolmogorov-Tihomirov metric dimension, and is a dimension in the 
 context of noncommutative geometric measure theory \cite{GI5,GI4}.
 \end{rem}

 \section{Asymptotic dimension of some cylindrical ends}
 \label{subsec:cylindrical}

 In this Section we want to compare our work with a work of 
 Davies.  In \cite{Davies} he defines the asymptotic dimension 
 of cylindrical ends of a Riemannian manifold $M$ as follows.  Let 
 $E\subset M$ be homeomorphic to $(1,\infty)\times A$, where $A$ is a 
 compact Riemannian manifold.  Set $\dE:=\{1\}\times A$, $E_r:=\{ x\in 
 E: \d(x,\dE) < r\}$, where $\d$ is the restriction of the metric in 
 $M$.  Then $E$ has asymptotic dimension $D$ if there is a positive 
 constant $c$ s.t.
 \begin{equation}\label{e:davies}
	 c^{-1} r^D \leq vol(E_r) \leq cr^D,
 \end{equation}
 for all $r\geq1$.  Davies does not assume bounded geometry for $E$.  If 
 one does, the two definitions coincide, more precisely if an 
 asymptotic dimension {\it \`a la} Davies exists, it coincides with 
 ours. 

 \begin{Prop}\label{2.3.1} 
	 With the above notation, if the volume form on $E$ is a uniformly 
	 bounded measure (as in Definition \ref{1.2.12}), or in particular 
	 if $E$ has bounded geometry (as in Definition \ref{2.1.1}), and 
	 there is $D$ as in {\rm (\ref{e:davies})}, then $\ad{E}=D$.
 \end{Prop}
 \begin{proof}
	 Choose $o\in E$, and set $\d:=\d(o,\dE)$, $\D:=diam(\dE)$.  Then 
	 it is easy to prove that $E_{R-\d-\D} \subset B_E(o,R) \subset 
	 E_{R+\d}$.  \\
	 Then $c^{-1}(R-\d-\D)^D \leq vol(B_E(o,R)) \leq c(R+\d)^D$, and 
	 from Proposition \ref{1.2.14} the thesis follows.  
 \end{proof}

 Motivated by (\cite{Davies}, example 16), let us set the following

 \begin{Dfn}\label{2.3.2} 
	 $E$ is a standard end of local dimension $N$ if it is 
	 homeomorphic to $(1,\infty)\times A$, endowed with the metric 
	 $ds^2=dx^2+f(x)^2d\om^2$, and with the volume form 
	 $dvol=f(x)^{N-1}dxdvol_{\om}$, where $(A,\om)$ is an $(N-1)$-dimensional 
	 compact Riemannian manifold, and $f$ is an increasing smooth 
	 function.
 \end{Dfn}

 \begin{Prop}\label{2.3.3} 
	 The volume form on a standard end $E$ is a uniformly bound\-ed 
	 measure.  Therefore, if $E$ satisfies equation {\rm (\ref{e:davies})}, we 
	 get $\ad{E}=D$.
 \end{Prop}
 \begin{proof}
	 It is easy to show that, for $(x_0,p_0)\in E$, $r<x_0-1$, \\
	 \begin{align*}
		 [x_0-r/2,x_0+r/2]&\times B_A\left(p_0, 
		 \frac{r/2}{f(x_0+r/2)}\right) \subset B_E((x_0,p_0),r) \\
		 & \subset [x_0-r,x_0+r]\times B_A\left(p_0, 
		 \frac{r}{f(x_0-r)}\right)
	 \end{align*}
	 So that, with $V_{X}(x,r) := vol(B_{X}(x,r))$, 
	 \begin{align*}
		 \int_{x_0-r/2}^{x_0+r/2}f(x)^{N-1}dx\ & V_A\left(p_0, 
		 \frac{r/2}{f(x_0+r/2)}\right) \leq V_E((x_0,p_0),r) \\
		 &\leq \int_{x_0-r}^{x_0+r}f(x)^{N-1}dx\ V_A\left(p_0, 
		 \frac{r}{f(x_0-r)}\right)
	 \end{align*}
	 which implies
	 \begin{align*}
		 rf(x_0-r/2)^{N-1}\ V_A\left(p_0,\frac{r/2}{f(x_0+r/2)}\right) 
		 &\leq V_E((x_0,p_0),r) \\
		 &\leq 2rf(x_0+r)^{N-1}\ V_A\left(p_0,\frac{r}{f(x_0-r)}\right)
	 \end{align*}
	 As for $x_0\to\infty$, $V_A(p_0,\frac{r}{ f(x_0-r)})\sim c 
	 \left(\frac{r}{f(x_0-r)}\right)^{N-1}$, and the same holds for 
	 $V_A(p_0,\frac{r/2}{f(x_0+r/2)})$, we get the thesis.  
 \end{proof}

 \begin{Cor}\label{2.3.4}  
	 Let $E$ be the standard end given by $E := (1,\infty)\times 
	 S^{N-1}$, endowed with the metric 
	 $ds^2=dr^2+r^{2(D-1)/(N-1)}d\om^2$, and with the volume form 
	 $dvol=r^{D-1}drd^{N-1}\om$ $(${\rm \cite{Davies}}, example $16)$.  
	 Then $\ad{E}=D$.
 \end{Cor}

 \begin{rem}
	Observe that $\ad{M}$ makes sense for any metric space, hence for 
	any cylindrical end, while Davies' asymptotic dimension does not.  
	Indeed let $E := (1,\infty)\times S^1$, endowed with the metric 
	$ds^2=dr^2+f(r)^2d\om^2$, and with the volume form 
	$dvol=f(r)drd\om$, where $f(r):=\frac{d}{ dr}(r^2\log r)$.  Then 
	$\ad{E}=2$, but $vol(E_r)$ does not satisfy one of the 
	inequalities in (\ref{e:davies}).
 \end{rem}
 
 Before closing this section we observe that the notion of standard 
 end allows us to construct an example which shows that we could 
 obtain quite different results if we used $\liminf$ instead of 
 $\limsup$ in the definition of the asymptotic dimension.  It makes 
 use of the following function
 $$
 f(x)=
 \begin{cases}
	 \sqrt{x} &  x\in[1,a_1] \\
	 2+ b_{n-1}+ c_{n-1} + (x-a_{2n-1}) & x\in[a_{2n-1},a_{2n}] \\
	 2+ b_{n-1}+ c_{n} + \sqrt{ x-a_{2n}+1 } & x\in[a_{2n},a_{2n+1}]
 \end{cases}
 $$
 where $a_0:=0$, $a_{n}-a_{n-1}:=2^{2^n}$, $b_{n}:= \sum_{k=1}^{n} 
 \sqrt{ 2^{2^{2k+1}}+1 }$, $c_{n}:= \sum_{k=1}^{n} ( 2^{2^{2k}}-1 )$, 
 $n\geq1$.

 \begin{Prop}\label{2.3.5} 
	 Let $M$ be the Riemannian manifold obtained as a $C^\infty$ 
	 regularization of $C\cup_\f E$, where $C:= \{ (x,y,z)\in\br^3 : 
	 (x-1)^2+y^2+z^2=1,\ x\leq1 \}$, with the Euclidean metric, 
	 $E:=[1,\infty)\times S^1$, endowed with the metric 
	 $ds^2=dx^2+f(x)^2d\om^2$, and with the volume form 
	 $dvol=f(x)dxd\om$, where $\f$ is the identification of $\{ 
	 y^2+z^2=1,\ x=1 \}$ with $\{1\}\times S^1$.  Then the volume form 
	 is a uniformly bounded measure, $\ad{M}\geq2$ but 
	 $\underline{d}_\infty(M)\leq3/2$, where 
	 $\underline{d}_\infty(M):=\limr\liminf_{R\to\infty}\frac{\log 
	 n_r(B_M(x,R)) }{ \log R}$.
 \end{Prop}
 \begin{proof}
	 Set $o:=(0,0,0)\in M$, then it is easy to see that, for 
	 $n\to\infty$, $a_{n}\sim 2^{2^n}$, $b_{n}\sim c_{n} \sim 
	 2^{2^{2n}}$, and
	 \begin{align*}
		 area(B_M(o,a_{2n})) & \sim \frac12 a_{2n}^2 \\
		 area(B_M(o,a_{2n-1}) )& \sim \frac53 a_{2n-1}^{3/2} 
	 \end{align*}
	 so that, calculating the limit of $\frac{\log 
	 area(B_M(o,R))}{\log R}$ on the sequence $R=a_{2n}$ we get $2$, 
	 while on the sequence $R=a_{2n-1}$ we get $3/2$.  The thesis 
	 follows easily, using Proposition \ref{1.2.14}.  
 \end{proof}
 
 \section{The asymptotic dimension and the 0-th Novikov Shubin 
 invariant} \label{subsec:NSinvariant}

 In this Section we show that, for a class of open manifolds of 
 bounded geometry, the asymptotic dimension coincides with the $0$-th 
 Novikov-Shubin invariant defined in \cite{GI4}.  In all this Section 
 $M$ denotes a manifold of C$^{\infty}$-bounded geometry, $i.e.$ $M$ 
 has positive injectivity radius, and the curvature tensor is bounded 
 together with all its covariant derivatives.  We assume moreover that 
 $M$ satisfies
 
 \begin{Assump}\label{ass:CG}
	there are $A,\ C,\ C'>0$ s.t. for all $x\in M$, $r>0$,
	 \begin{align}
		 V(x,2r) & \leq A V(x,r) \label{e:voldouble}\\
		 \frac{C}{ V(x,\sqrt{r})} & \leq p_{r}(x,x) \leq 
		 \frac{C'}{V(x,\sqrt{r})}   \label{e:stime}
	 \end{align}
	 where $V(x,r) := vol(B(x,r))$ and $p_{t}(x,y)$ is the heat kernel on $M$.
 \end{Assump}
 
 \begin{rem}
	  $(i)$ Inequality~(\ref{e:voldouble}) is introduced in 
	  \cite{CG96} and called the volume doubling property.  \\
	  $(ii)$ A result of Coulhon-Grigor'yan (\cite{CG96}, Corollary 
	  7.3) (\cite{Grigoryan94}, Proposition 5.2) states that 
	  Assumption \ref{ass:CG} is equivalent to the following 
	  isoperimetric inequality introduced in \cite{CG96}.  There are 
	  $\a,\ \b>0$ s.t. for all $x\in M$, $r>0$, and all regions 
	  $U\subset B(x,r)$,
	  $$
	  \l_1(U) \geq \frac{\a }{ r^2} \left( \frac{V(x,r)}{ vol(U)} 
	  \right)^\b,
	  $$
	  where $\l_1(U)$ is the first Dirichlet eigenvalue of $\D$ in 
	  $U$.  \\
	  $(iii)$ Assumption \ref{ass:CG} is satisfied by all manifolds 
	  with positive Ricci curvature \cite{LY}, and covering manifolds 
	  whose group of deck transformations has polynomial growth \cite{SC}. 
 \end{rem}

 \begin{Prop}\label{2.1.7} 
	 Let $M$ be a complete Riemannian manifold of 
	 C$^{\infty}$-bounded geometry, satisfying Assumption 
	 \ref{ass:CG}.  \\
	 Then $\ad{M}=\limsup_{t\to\infty} \frac{-2\log 
	 p_{t}(x,x)}{\log t}$, for any $x\in M$.
 \end{Prop}
 \begin{proof} 
	 Follows from Proposition \ref{2.1.3} and estimates (\ref{e:stime}).
 \end{proof}
 
 \begin{rem}
	The previous result shows that there are some connections between 
	the asymptotic dimension of a manifold and the notion of dimension 
	at infinity for semigroups (in our case the heat kernel semigroup) 
	considered by Varopoulos (see \cite{VSC}).
 \end{rem}

 The volume doubling property is a weak form of polynomial growth 
 condition, and still guarantees the finiteness of the asymptotic 
 dimension (for manifolds of bounded geometry).

 \begin{Prop}
	 Let $M$ be a complete Riemannian manifold of 
	 C$^{\infty}$-bounded geometry, and suppose the volume 
	 doubling property (\ref{e:voldouble}) holds.  Then $M$ has 
	 finite asymptotic dimension.
 \end{Prop}
 \begin{proof}
	 Let $R>1$, and $n\in\bn$ be s.t. $2^{n-1}<R \leq 2^n$.  Then 
	 $V(x,R)\leq V(x,2^n)\leq A^n V(x,1)$, so that
	 $$
	 1 \leq \frac{V(x,R) }{ V(x,1)} \leq A^n \leq A R^{\log_2 A}.
	 $$
	 Therefore $\ad{M} = \limsup_{R\to\infty} \frac{\log V(x,R) }{ 
	 \log R} \leq \log_2 A$.  
 \end{proof}

 \begin{Dfn} {\rm \cite{Roe1}}
	A regular exhaustion $\ck$ is an increasing sequence $\{K_{n}\}$ 
	of compact subsets of $M$, whose union is $M$, and such that, for 
	any $r>0$
	$$									
	\lim_{n\to\infty} \frac{vol(Pen^{+}(K_{n},r))} 
	{vol(Pen^{-}(K_{n},r))} =1,
	$$
	where we set $Pen^{+}(K,r):= \{x\in M: \d(x,K)\leq r\}$, and 
	$Pen^{-}(K,r):=$ the closure of $M\setminus Pen^{+}(M\setminus 
	K,r)$.
 \end{Dfn}

 \begin{Prop}\label{p:exbyballs} 
	 Let $M$ be an open manifolds of \bg and satisfying Assumption 
	 \ref{ass:CG}. 
	 \itm{i} There is $\g>0$ s.t. for any $x,\ y\in M$, $r>0$, if 
	 $B(x,r)\cap B(y,r) \neq \emptyset$, then
	 $$
	 \g^{-1}\leq\frac{V(x,r)}{V(y,r)}\leq\g.
	 $$
	 \itm{ii} There is a sequence $n_{k}\in\bn$ s.t. $\{B(x,n_{k})\}$ 
	 is a regular exhaustion of $M$.
 \end{Prop}
 \begin{proof} 
	 $(i)$ The inequality easily follows by a result of Grigor'yan 
	 (\cite{Grigoryan94}, Proposition 5.2), where it is shown that 
	 Assumption \ref{ass:CG} implies the existence of a 
	 constant $\g$ such that
	 $$
	 \g^{-1}\left(\frac{R}{r}\right)^{\a_{1}} \leq 
	 \frac{V(x,R)}{V(y,r)} \leq \g\left(\frac{R}{r}\right)^{\a_{2}}
	 $$
	 for some positive constants $\a_{1},\ \a_{2}$, for any $R\geq r$, 
	 and $B(x,R)\cap B(y,r)\neq\emptyset$. 
	 \itm{ii} The statement follows from the fact that the volume 
	 doubling property implies subexponential (volume) growth, so that 
	 the result is contained in (\cite{Roe1}, Proposition 6.2).
 \end{proof}
  
 Recall from \cite{GI3} that the C$^{*}$-algebra $\ca$ of almost 
 local operators on $M$ is the norm closure of the finite propagation 
 operators on $L^{2}(M,dvol)$. Then
 
 \begin{Prop}\label{thm:trace} {\rm \cite{GI3}}
	 There is on $\ca$ a lower semicontinuous semifinite trace 
	 $Tr_{\ck}$, which, on the heat semigroup, is given by the 
	 following formula,
	 $$
	 Tr_{\ck}(\e{-t\D}) = 
	 \Lim_{\om}\frac{\int_{K_{n}}tr(p_{t}(x,x))dvol(x)}{vol(K_{n})},
	 $$
	 where $\Lim_{\om}$ is a generalized limit.
 \end{Prop}
 
 \begin{rem}
     $(i)$  The above formula for the trace was considered by J. Roe 
     in \cite{Roe1}.  However, this formula does not describe a 
     semicontinuous trace on the C$^*$-algebra of almost local 
     operators.  Therefore we introduced a semicontinuous semifinite 
     regularization in \cite{GI3}.
     \\
     $(ii)$ $L^{2}$-Betti numbers for open manifolds have been 
     introduced in \cite{Roe1}, where it is shown that the $0$-th 
     $L^{2}$-Betti number of a noncompact manifold is zero.  For this 
     reason it does not appear in the formula for $\a_{0}$ below.
 \end{rem}

 By means of $Tr_{\ck}$ we defined the $0$-th Novikov-Shubin invariant 
 as
 $$
 \a_{0}(M,\ck) := 2\limsup_{t\to\infty} \frac{\log Tr_\ck(\e{-t\D})}
 {\log 1/t}.
 $$

 \begin{Thm}\label{a0=dinfinity} 
	 Let $M$ be an open manifold of \bg and satisfying Assumption 
	 \ref{ass:CG}, endowed with the regular 
	 exhaustion $\ck$ given by Proposition \ref{p:exbyballs} $(ii)$.  
	 Then the asymptotic dimension of $M$ coincides with the 0-th 
	 Novikov-Shubin invariant, namely $d_\infty(M)=\a_0(M,\ck)$. In 
	 particular $\a_{0}$ is independent of the limit procedure $\Lim_{\om}$.
 \end{Thm}
 \begin{proof} 
	 First, from equation (\ref{e:stime}) and Proposition 
	 \ref{p:exbyballs} $(i)$, we get
	 \begin{align*}
		 \frac{C\g^{-1}}{V(o,\sqrt{t})} &\leq \frac {\int_{B(o,r)} 
		 \frac{C}{V(x,\sqrt{t})} dvol(x)} {V(o,r)} \leq 
		 \frac{\int_{B(o,r)}p_{t}(x,x)dvol(x)}{V(o,r)}\\
		 &\leq \frac{\int_{B(o,r)} \frac{C'}{V(x,\sqrt{t})} dvol(x)} 
		 {V(o,r)} \leq \frac{C'\g}{V(o,\sqrt{t})}
	 \end{align*}
	 therefore, by Proposition \ref{thm:trace} we have,
	 $$
	 \frac{C\g^{-1}}{V(o,\sqrt{t})} \leq Tr_{\ck}(\e{-t\D}) \leq 
	 \frac{C'\g}{V(o,\sqrt{t})}
	 $$
	 hence, finally,
	 \begin{align*}
		 d_\infty(M) & = 
		 2\limsup_{t\to\infty}\frac{\log(V(o,t))}{2\log t} = 
		 2\limsup_{t\to\infty}\frac{\log(C'\g
		 V(o,\sqrt{t})^{-1})}{\log\frac{1}{t}}\\
		 \leq \a_0(M,\ck) & \equiv 2\limsup_{t\to\infty}\frac 
		 {\log Tr_{\ck}(\e{-t\D})}{\log\frac{1}{t}} 
		 \leq2\limsup_{t\to\infty}\frac {\log(C\g^{-1}  
		 V(o,\sqrt{t})^{-1})}{\log\frac{1}{t}}\\
		 & = 2\limsup_{t\to\infty}\frac {\log(V(o,t))}{2\log t} = 
		 d_\infty(M).
	 \end{align*}
	 The thesis follows.
 \end{proof}

 \begin{ack}
	 We would like to thank D.~Burghelea, M.~Farber, L.~Friedlander 
	 for conversations.  We also thank I.~Chavel, E.B.~Davies, 
	 A.~Grigor'yan, P.~Li, L.~Saloff-Coste for having suggested useful 
	 references on heat kernel estimates, and the referee for drawing 
	 our attention to the notion of Gromov asymptotic dimension. 
 \end{ack}



\begin{thebibliography}{99}
     
 \bibitem{Atiyah} M. F. Atiyah.  {\it Elliptic operators, discrete 
 groups and von Neumann algebras}.  Soc.  Math.  de France, 
 Ast\'erisque {\bf 32--33} (1976), 43--72.
 
 \bibitem{Chavel} I. Chavel.  {\it Riemannian geometry - A modern 
 introduction}.  Cambridge Univ.  Press, Cambridge, 1993.

 \bibitem{Co} A. Connes.  {\it Non Commutative Geometry}.  Academic 
 Press, 1994.

 \bibitem{CG96} T. Coulhon, A. Grigor'yan.  {\it On-diagonal lower 
 bounds for heat kernels and Markov chains}.  Duke Math. J., {\bf 89} 
 (1997), 133--199.

 \bibitem{Davies} E. B. Davies.  {\it Non-gaussian aspects of heat 
 kernel behaviour}.  J. London Math. Soc., {\bf 55} (1997), 105--125.

 \bibitem{Dra} A.N. Dranishnikov.  {\it Asymptotic topology}.  
 Preprint math/9907192.
 
 \bibitem{Grigoryan94} A. Grigor'yan.  {\it Heat kernel upper bounds 
 on a complete non-compact manifold}.  Revista Matematica 
 Iberoamericana, {\bf 10} (1994), 395--452.

 \bibitem{Gromov} M. Gromov.  {\it Asymptotic invariants of infinite 
 groups}.  Geometric group theory, Vol.  2 (Sussex, 1991), 1--295, 
 London Math.  Soc.  Lecture Note Ser., 182, Cambridge Univ.  Press, 
 Cambridge, 1993.

 \bibitem{GS} M. Gromov, M. Shubin.  {\it Von Neumann spectra near 
 zero}.  Geometric and Functional Analysis, {\bf 1} (1991), 375--404.

 \bibitem{GI1} D. Guido, T. Isola.  {\it Singular traces for 
 semifinite von~Neumann algebras}.  Journal of Functional Analysis, 
 {\bf 134} (1995), 451--485.
 
 \bibitem{GI5} D. Guido, T. Isola.  {\it Singular traces, dimensions, 
 and Novikov-Shubin invariants}.  Operator Theoretical Methods, 
 Proceedings of the 17th Conference on Operator Theory, (Timisoara, 
 Romania, 1998) 151--171, A. Gheondea, R.N. Gologan, D. Timotin 
 Editors, The Theta Foundation, Bucharest 2000.

 \bibitem{GI3} D. Guido, T. Isola.  {\it A semicontinuous trace for 
 almost local operators on an open manifold}.  Preprint.

 \bibitem{GI4} D. Guido, T. Isola.  {\it Noncommutative Riemann 
 integration and singular traces for C$^{*}$-algebras}.  Journ.  
 Funct.  Analysis, {\bf 176} (2000), 115-152.

 \bibitem{KT} I. Kolmogoroff, M.G. Tihomirov.  {\it $\eps$-Entropy and 
 $\eps$-capacity of sets in functional spaces}.  A. M. S. Transl., 
 {\bf 17} (1961), 277--364.
 
 \bibitem{LY} P. Li, S. T. Yau. {\it On the parabolic kernel of the 
 Schr\"odinger operator}. Acta Math., {\bf 156} (1986), 153--201.

 \bibitem{NS} S. P. Novikov, M. A. Shubin.  {\it Morse theory and von 
 Neumann {\rm II}${}_1$ factors}.  Doklady Akad.  Nauk SSSR, {\bf 289} 
 (1986), 289--292.\\
 S. P. Novikov, M. A. Shubin.  {\it Morse theory and von 
 Neumann invariants on non-simply connected manifolds}.  Uspekhi Math.  
 Nauk, {\bf 41}, 5 (1986), 222--223 (in Russian).

 \bibitem{Pontriagin} A. V. Arkhangel'skii, L. S. Pontryagin (eds.).  
 {\it General Topology I. Basic concepts and constructions.  Dimension 
 theory}.  Enc.  Math.  Sci.  {\bf 17}, Springer, New York, 1988.
 
 \bibitem{Roe1} J. Roe.  {\it An index theorem on open manifolds.  I, 
 II}.  J. Diff.  Geom., {\bf 27} (1988), 87--136.

 \bibitem{Salli} A. Salli.  {\it On the Minkowski dimension of 
 strongly porous fractal sets in $\br^{n}$}.  Proc.  London Math.  
 Soc., {\bf 62} (1991), 353--372.
 
 \bibitem{SC} L. Saloff-Coste.  {\it A note on Poincar\'e, Sobolev, 
 and Harnack inequalities}.  Internat.  Math.  Res.  Notices (1992), 
 no.  2, 27--38.

 \bibitem{Varopoulos1} N. T. Varopoulos.  {\it Random walks and 
 Brownian motion on manifolds}.  Symposia Mathematica, {\bf XXIX} 
 (1987), 97--109.

 \bibitem{VSC} N. T. Varopoulos, L. Saloff-Coste, T. Couhlon.  {\it 
 Analysis and geometry on groups}.  Cambridge Univ.  Press, Cambridge, 
 1992.
 
\bibitem{Yu} Guoliang Yu. {\it The Novikov conjecture for groups with 
finite asymptotic dimension}.  Ann.  of Math.  (2) {\bf 147} (1998), 
325--355.

\end{thebibliography}
\end{document}